\documentclass[10pt]{amsart}
\usepackage{amsthm,amsmath,amssymb,amsfonts}
\usepackage[american]{babel}
\usepackage{url}

\newtheorem{theorem}{Theorem}
\newtheorem*{theoremA}{Theorem~A}
\newtheorem*{theoremB}{Theorem~B}
\newtheorem{lemma}[theorem]{Lemma}

\newtheorem{corollary}[theorem]{Corollary}

\theoremstyle{definition}

\newtheorem{example}[theorem]{Example}

\theoremstyle{remark}

\begin{document}


\newcommand{\NU}[1]{\mbox{\rm V}(#1)}
\newcommand{\lara}[1]{\langle{#1}\rangle}
\newcommand{\ZZ}{\mbox{$\mathbb{Z}$}}
\newcommand{\sZZ}{\mbox{$\scriptstyle\mathbb{Z}$}}
\newcommand{\QQ}{\mbox{$\mathbb{Q}$}}
\newcommand{\paug}[2]{\varepsilon_{#1}(#2)}


\title[Unit groups with no noncyclic abelian finite subgroups]
{Unit groups of integral finite group rings
with no noncyclic abelian finite subgroups}

\author{Martin Hertweck}
\address{Universit\"at Stuttgart, Fachbereich Mathematik, IGT,
Pfaffenwaldring 57, 70550 Stuttgart, Germany}
\email{hertweck@mathematik.uni-stuttgart.de}



\date{\today}

\begin{abstract}
It is shown that in the units of augmentation one of an
integral group ring $\ZZ G$ of a finite group $G$, a noncyclic 
subgroup of order $p^{2}$, for some odd prime $p$, exists only 
if such a subgroup exists in $G$.
The corresponding statement for $p=2$ holds by the Brauer--Suzuki 
theorem, as recently observed by W.~Kimmerle.
\end{abstract}

\maketitle

\section{Introduction}\label{Sec:intro}

Is a finite subgroup $H$ of units in the integral group ring $\ZZ G$ of 
a finite group $G$ necessarily isomorphic to a subgroup of $G$?
Of course, torsion coming from the coefficient ring should be excluded,
that is, only finite subgroups $H$ in $\NU{\ZZ G}$, the group of units of 
augmentation one in $\ZZ G$, will be considered.
The question was raised by Higman in his thesis (1940),
where he gave an affirmative answer when $G$ is metabelian nilpotent
or the affine group over a prime field; cf.\ Sandling (1981).
In the survey of Sandling (1984) it is included as Problem~5.4,
and noted that an affirmative answer for metabelian 
$G$ was finally given by Roggenkamp (1981); but see also 
Cliff, Sehgal and Weiss (1981), and Marciniak and Sehgal (2003) 
for a more recent result, giving a generalization based on 
a theorem of Weiss (1988).
These results are really about certain `large' 
torsion-free normal subgroups of $\NU{\ZZ G}$. For a more complete
discussion, see Chapter~4 in Sehgal's book (1993).

As a sort of converse, one may fix a finite group $H$ and look for
groups $G$ for which $H$ embeds into $\NU{\ZZ G}$, again hoping for the 
best, but little is known in this respect.
What is known is that if a cyclic group $H$ of prime power order
embeds into some unit group $\NU{\ZZ G}$, then $H$ also embeds into $G$
(due to an observation of Cohn and Livingstone (1965); see also
Zassenhaus (1974)), and only recently in Hertweck (2007b)
it was shown that the restriction on the order can be removed if in 
addition $G$ is assumed to be solvable.
In this spirit, Marciniak, at a satellite conference of the ICM 2006,
asked whether a group $G$ necessarily has a subgroup isomorphic to 
Klein's four group provided this is the case for $\NU{\ZZ G}$.
Kimmerle immediately observed that this is implied by the Brauer--Suzuki
theorem (rendered in Kimmerle (2006)), see Section~\ref{Sec:Kim}.
Our complementary result is as follows.
\begin{theoremA}
Let $G$ be a finite group. Suppose that $\NU{\ZZ G}$ has a noncyclic
abelian subgroup of order $p^{2}$, for some odd prime $p$. Then the same 
is true for $G$ (i.e., Sylow $p$-subgroups of $G$ are not cyclic).
\end{theoremA}

It is easy to verify that a finite $p$-group with no noncyclic abelian
subgroup is either cyclic or a (generalized) quaternion group,
see Theorem~4.10 in Gorenstein (1968).
It comes to mind that the theory of cyclic blocks might be used in the
proof, but it is pretty simple and makes only use of a fact about 
vanishing of partial augmentations of torsion units, established in
Hertweck (2006, 2007a).

We remark that both results (whether $p$ is even or odd) for a solvable 
group $G$ are covered by Theorem~5.1 in Dokuchaev and Juriaans (1996).

Note that a group $G$ whose Sylow $2$-subgroups are cyclic has a normal
$2$-comple\-ment, by Burnside's well known criterion, see Theorem~4.3 in
Gorenstein (1968). We obtain the following corollary.

\begin{corollary}\label{C1}
Let $G$ be a finite group having cyclic Sylow $p$-subgroups 
for some prime $p$. Then any finite $p$-subgroup of $\NU{\ZZ G}$
is isomorphic to a subgroup of $G$.
\end{corollary}

Finally, we remark that, as with other results in this field, the 
theorem can be formulated for more general coefficient rings than
$\ZZ$, notably for the semilocalization of $\ZZ$ at the prime divisors
of the order of $G$. Unfortunately, it is 
definitely wrong for $p$-adic coefficient rings.

\section{Kimmerle's observation}\label{Sec:Kim}

Coming back to the initial question, we mention that in the hope for 
further positive results, it is natural to impose restrictions on the 
prime divisors of the finite subgroup $H$, i.e., to consider only 
$\pi$-groups $H$ for some set $\pi$ of primes (a singleton $\{p\}$, to
begin with), as has been done before in work on the stronger Zassenhaus
conjecture (ZC3), cf.\ Dokuchaev and Juriaans (1996). It is well 
known that then, one can assume that $\text{\rm O}_{\pi^{\prime}}(G)$,
the largest normal $\pi^{\prime}$-subgroup of $G$, is trivial, for $H$
has an isomorphic image under the natural map 
$\ZZ G\rightarrow\ZZ G/\text{\rm O}_{\pi^{\prime}}(G)$,
see the remark after Theorem~2.2 in Dokuchaev and Juriaans (1996).

This derives from the vanishing of certain partial augmentations of
the elements of $H$. Recall that for a group ring element
$u=\sum_{g\in G}a_{g}g$ (all $a_{g}$ in $\ZZ$), its partial 
augmentation with respect to an element $x$ of $G$, or rather its
conjugacy class $x^{G}$ in $G$, is the sum $\sum_{g\in x^{G}}a_{g}$;
we will denote it by $\paug{x}{u}$. The result of Cohn and Livingstone 
mentioned in the introduction really says that if an element $h$ of $H$
is of prime power order, then there exists an element $x$ in $G$ of
the same order such that $\paug{x}{h}\neq 0$.
Note that $\paug{z}{u}=a_{z}$ for an element $z$ in the center of $G$.
An old yet fundamental result from Berman (1955) and Higman (1940)
asserts that if $\paug{z}{h}\neq 0$ for an element $h$ in $H$ and some 
$z$ in the center of $G$, then $h=z$.

Coming to Marciniak's question, suppose that $G$ has no subgroups 
isomorphic to Klein's four group. For our purpose, we can assume that
$\text{\rm O}_{2^{\prime}}(G)=1$ and that Sylow $2$-subgroups of $G$
are not cyclic. Thus Sylow $2$-subgroups of $G$ are (generalized)
quaternion, and by the Brauer--Suzuki theorem, from Brauer and Suzuki
(1959), $G$ contains a unique involution $z$. For an involution $u$ in
$\NU{\ZZ G}$, the Cohn--Livingstone result gives $\paug{z}{u}\neq 0$,
and therefore $u=z$ by the Berman--Higman result, answering 
Marciniak's question in the affirmative.

\begin{theoremB}[Kimmerle]
Let $G$ be a finite group. Suppose that $\NU{\ZZ G}$ has a subgroup 
isomorphic to Klein's four group. Then the same is true for $G$. 
\end{theoremB}

We do not know of a proof avoiding the use of the Brauer--Suzuki theorem.

Suppose that Sylow $2$-subgroups of $G$ are quaternion groups. Then
the theorem implies that finite $2$-subgroups of $\NU{\ZZ G}$ are
cyclic or quaternion groups. Taking into account the structure of the 
quaternion groups, and the Cohn--Livingstone result, one obtains
the following corollary.

\begin{corollary}\label{C2}
Let $G$ be a finite group whose Sylow $2$-subgroups are quaternion groups
(ordinary or generalized). Then any finite $2$-subgroup of $\NU{\ZZ G}$
is isomorphic to a subgroup of $G$.
\end{corollary}

\section{Proof of Theorem~A}\label{Sec:proof}

The partial augmentations of a torsion unit in $\NU{\ZZ G}$ encode
its character values in a way establishing a connection to group
elements which respects a divisibility relation between orders.
We will make use of a lemma which is an easy consequence of this fact.

\begin{lemma}\label{L1}
Let $u$ be a torsion unit in $\NU{\ZZ G}$ of, say, order $n$. Let $s$ be
a natural integer coprime to $n$, so that $st\equiv 1\mod{n}$ for another
natural integer $t$. Then for all $x$ in $G$ whose order divide $n$, we
have $\paug{x}{u^{s}}=\paug{x^{t}}{u}$.
\end{lemma}
\begin{proof}
Let $\zeta$ be a primitive $n$-th complex root of unity, and let $\sigma$
be the Galois automorphism of $\QQ(\zeta)$ sending $\zeta$ to $\zeta^{s}$.
Let $x_{1},\dotsc,x_{k}$ be representatives of the conjugacy classes
of $G$ whose elements have order dividing $n$. Note that then
$x_{1}^{t},\dotsc,x_{k}^{t}$ is another system of representatives. By 
Theorem~2.3 in Hertweck (2007a), $\paug{x}{u}\neq 0$ is possible only for 
elements $x$ whose order divide $n$. Thus for any ordinary irreducible 
character $\chi$ of $G$, we have
\begin{align*}
& \sum_{i=1}^{k}\paug{x_{i}}{u^{s}}\chi(x_{i})
= \chi(u^{s}) = \chi(u)^{\sigma} 
= \sum_{i=1}^{k}\paug{x_{i}}{u}\chi(x_{i})^{\sigma}
\\ & \qquad
= \sum_{i=1}^{k}\paug{x_{i}}{u}\chi(x_{i}^{s})
= \sum_{i=1}^{k}\paug{x_{i}^{t}}{u}\chi(x_{i}).
\end{align*}
Since the character table of $G$, stripped off from any additional 
information, is an invertible matrix, it follows that
$\paug{x_{i}}{u^{s}}=\paug{x_{i}^{t}}{u}$ for all indices $i$, which 
proves the lemma.
\end{proof}

\begin{corollary}\label{C4}
Let $u$ be a torsion unit in $\NU{\ZZ G}$ of, say, order $n$. Then
for any $x$ in $G$ whose order divides $n$,
\[ \sum_{s\in(\sZZ/n\sZZ)^{\times}}\paug{x}{u^{s}}
= \sum_{s\in(\sZZ/n\sZZ)^{\times}}\paug{x^{s}}{u}. \]
\end{corollary}

\begin{corollary}\label{C5}
Suppose that for a prime divisor $p$ of the order of $G$, all elements of
order $p$ in $G$ are conjugate to a power of some fixed element $x$.
Let $u$ be a torsion unit in $\NU{\ZZ G}$ of order $p$. Then
$\sum_{i=1}^{p-1}u^{i}$ and $\sum_{i=1}^{p-1}x^{i}$ have the same partial 
augmentations.
\end{corollary}
\begin{proof}
Let $k$ be the number of conjugacy classes of elements of order $p$ in $G$.
By Corollary~\ref{C4} and Theorem~2.3 in Hertweck (2007a),
\[ \varepsilon_{x}\Bigg(\sum_{i=1}^{p-1}u^{i}\Bigg) = 
\sum_{i=1}^{p-1}\paug{x^{i}}{u} = 
\frac{p-1}{k}\sum_{y^{G}\colon y\in\lara{x}}\paug{y}{u} =
\frac{p-1}{k} = \varepsilon_{x}\Bigg(\sum_{i=1}^{p-1}x^{i}\Bigg). \]
Applying again Theorem~2.3 from Hertweck (2007a), the corollary follows.
\end{proof}

We will apply this by means of the following formula
relating ranks of an idempotent to arithmetical properties of the group.

\begin{corollary}\label{C6}
Suppose that for a prime divisor $p$ of the order of $G$, all elements of
order $p$ in $G$ are conjugate to a power of some fixed element $x$.
Suppose further that $\NU{\ZZ G}$ contains an elementary abelian subgroup
$U$ of order $p^{2}$. Then for any ordinary character $\chi$ of $G$,
\begin{equation}\label{E1} 
\chi\Bigg(\frac{1}{p^{2}}\sum_{u\in U}u\Bigg) = 
\frac{1}{p^{2}}\Bigg(\chi(1) + (p+1) \sum_{i=1}^{p-1}\chi(x^{i})\Bigg). 
\end{equation}
\end{corollary}

We now turn to the proof of Theorem~A.
Suppose that $G$ has a cyclic Sylow $p$-subgroup $P$ ($p=2$ is allowed).
Let $x$ be an element of order $p$ in $P$, and set 
$N=\text{N}_{G}(\lara{x})$.
Suppose further that $\NU{\ZZ G}$ contains an elementary abelian subgroup
$U$ of order $p^{2}$.
Let $\chi$ be the character of $G$ which is induced from the principal
irreducible character of $P$. Then the rank in \eqref{E1} is
\[ \frac{1}{p^{2}}(|G:P| + |N:P|(p^{2}-1)). \]
If $\chi$ is a character of $G$ which is induced from a faithful
irreducible character of $P$, the rank in \eqref{E1} is
\[ \frac{1}{p^{2}}(|G:P| - |N:P|(p+1)). \]
The difference of these ranks is $|N:P|(p^{2}+p)/p^{2}$, which is 
not an integer. This contradiction proves the theorem.

In view of Corollaries~\ref{C1} and \ref{C2}, one may be tempted to
investigate the analogous problem for groups with dihedral Sylow 
$2$-subgroups. These groups were classified by Gorenstein and Walter,
and listed, for example, on p.~462 in Gorenstein (1968).
To indicate what can be done by now, we end with an example.

Note that the order of a finite subgroup of $\NU{\ZZ G}$ divides the
order of $G$, see Lemma~37.3 in Sehgal (1993); a fact which, surprisingly
enough from today's point of view, is in this
generality not recorded in Higman's thesis.

\begin{example}
For the alternating group $A_{7}$, any finite $2$-subgroup of 
$\NU{\ZZ A_{7}}$ is isomorphic to a subgroup of $A_{7}$.
\end{example}
\begin{proof}
Sylow $2$-subgroups of $A_{7}$ are dihedral of order $8$.
Let $x$ be an element of order $4$ in $A_{7}$. Then $x^{G}$ and 
$(x^{2})^{G}$ are the only conjugacy classes of elements of order $4$ 
and $2$, respectively. 
There is an (irreducible) character $\chi$ of $A_{7}$ of degree $6$
which is afforded by a deleted permutation representation. We have 
$\chi(x)=0$ and $\chi(x^{2})=2$. 

Let $U$ be a finite $2$-subgroup of $\NU{\ZZ A_{7}}$.
If $U$ is of order $2$, then $U$ is rationally conjugate to a subgroup
of $A_{7}$ by Corollary~3.5 in Hertweck (2006). If $U$ is of order
$4$, the Luthar--Passi method as described in Hertweck (2007a) is not 
sufficient to guarantee rational conjugacy to a subgroup of $A_{7}$:
for a unit $u$ of order $4$ in $\NU{\ZZ A_{7}}$ one cannot exclude the
possibility of having $(\paug{x^{2}}{u},\paug{x}{u})=(2,-1)$ when
$\chi(u)=4$. In this case, also $\chi(u^{-1})=4$. Anyway, $U$ is 
isomorphic to a subgroup of $A_{7}$, and the same is true
if $U$ is a Klein's four group.

Suppose that $U$ is abelian of order $8$. By the Cohn--Livingstone
result, $U$ is not cyclic. Set $e=\frac{1}{8}\sum_{u\in U}u$. Since $e$ 
is an idempotent, $\chi(u)$ is a rational integer. If $U$ is elementary
abelian, then $\chi(e)=\frac{1}{8}(\chi(1)+7\chi(x^{2}))=\frac{20}{8}$,
which is impossible. Thus $U$ contains $3$ elements of order $2$ and $4$
elements of order $4$. Trying out all possibilities shows that again 
$\chi(e)$ is not a rational integer.

It remains to consider the case when $U$ is the quaternion group. Let 
$u$ be an element of order $4$ in $U$. Since $\chi(u^{2})=\chi(x^{2})$,
the restriction of the character $\chi$ to $U$ is the sum of four linear
characters and the one of degree two. But this is impossible since 
$\chi$ is afforded by a rational representation, while the character
of degree two of the quaternion group comes from the block of the 
rational quaternion algebra (whence the name of the group).
\end{proof}



\section*{References}

{\small
\begin{list}{}{
\setlength{\leftmargin}{\parindent}
\setlength{\labelwidth}{\parindent}
\addtolength{\labelwidth}{-\labelsep}}
\item[Berman, S.~D. (1955).] 
On the equation {$x\sp m=1$} in an integral group ring.
{\it Ukrain. Mat. \v Z.} 7:253--261.
\item[Brauer, R., Suzuki, M. (1959).]
On finite groups of even order whose $2$-Sylow group is a quaternion group.
{\it Proc. Nat. Acad. Sci. U.S.A.} 45:1757--1759.
\item[Cliff, G.~H., Sehgal, S.~K., Weiss, A.~R. (1981).]
Units of integral group rings of met\-abel\-ian groups.
{\it J. Algebra} 73(1):167--185.
\item[Cohn, J.~A., Livingstone, D. (1965).]
On the structure of group algebras.~I.
{\it Canad. J. Math.} 17:583--593.
\item[Dokuchaev, M.~A., Juriaans, S.~O.~(1996).] 
Finite subgroups in integral group rings.
{\it Can\-ad. J. Math.} 48(6):1170--1179.
\item[Gorenstein, D. (1968).]
{\it Finite groups.} New York: Harper \& Row.
\item[Hertweck, M. (2006).] 
On the torsion units of some integral group rings.
{\it Algebra Colloq.} 13(2):329--348.
\item[Hertweck, M. (2007a).] 
Partial augmentations and Brauer character values of
torsion units in group rings.
{\it Comm.\ Algebra}, to appear (e-print \url{arXiv:math.RA/0612429v2}).
\item[Hertweck, M. (2007b).] 
The orders of torsion units in integral group rings of
finite solvable groups.
{\it Comm.\ Algebra}, to appear (e-print \url{arXiv:math.RT/0703541}).
\item[Higman, G. (1940).]
Units in group rings. {\it Ph.D. thesis.} University of Oxford
(Balliol College).
\item[Kimmerle, W. (2006).]
Arithmetical properties of finite groups. Talk
delivered at the Math Colloquium of the Vrije Universiteit Brussel.
\item[Marciniak, Z., Sehgal, S.~K. (2003).]
The unit group of $1+\Delta(G)\Delta(A)$ is torsion-free.
{\it J. Group Theory} 6(2):223--228.
\item[Roggenkamp, K.~W. (1981).]
Units in integral metabelian grouprings. I.
Jackson's unit theorem revisited.
{\it Quart. J. Math. Oxford Ser. (2)} 32:209--224.
\item[Sandling, R. (1981).]
Graham Higman's thesis ``Units in group rings''.
In: {\it Integral representations and applications (Oberwolfach, 1980).} 
Lecture Notes in Math.\ Vol.~882. Berlin: Springer, pp.~93--116.
\item[Sandling, R. (1984).]
The isomorphism problem for group rings: a survey.
In: {\it Orders and their applications (Oberwolfach, 1984).}
Lecture Notes in Math.\ Vol.~1142.
Berlin: Springer, pp.~256--288.
\item[Sehgal, S.~K. (1993).]
{\it Units in integral group rings}.
Pitman Monographs and Surveys in Pure and Applied Mathematics Vol.~69.
Harlow: Longman Scientific \& Technical.
\item[Weiss, A. (1988).]
Rigidity of {$p$}-adic {$p$}-torsion.
{\it Ann. of Math. (2)} 127(2):317--332.
\item[Zassenhaus, H. (1974).]
On the torsion units of finite group rings.
In: {\it Studies in mathematics (in honor of A. Almeida Costa).} 
Lisbon: Instituto de Alta Cultura, pp.~119--126.
\end{list}
}

\end{document}